\documentclass[11pt]{article}
\usepackage[dvips]{graphicx}
\usepackage{latexsym}
\usepackage{amssymb}
\newcommand{\C}{\mathbb{C}}
\newcommand{\CP}{\mathbb{CP}}

\newcommand{\R}{\mathbb{R}}
\newcommand{\Z}{\mathbb{Z}}
\newcommand{\G}{\mathbb{G}}

\newcommand{\koniec}{\begin{flushright}  $\Box $ \end{flushright}}
\def\be{\begin{equation}}
\def\ee{\end{equation}}

\def\Sm{\Sigma}

\def\Om{\Omega}
\def\Th{\Theta}

\def\om{\omega}

\def\p{\partial}

\newcommand{\hook}{{\setlength{\unitlength}{11pt}   
                   \begin{picture}(.833,.8)
                   \put(.15,.08){\line(1,0){.35}}
                   \put(.5,.08){\line(0,1){.5}}
                   \end{picture}}}
\def\a{\alpha}
\def\l{\lambda}

\def\tom{\tilde{\omega}}

\newtheorem{prop}{Proposition}  
\newtheorem{lemma}[prop]{Lemma}

\begin{document}

\title{The Twisted Photon Associated to Hyper-Hermitian
Four-Manifolds}

\author{Maciej Dunajski\thanks{email: dunajski@maths.ox.ac.uk}
\\ Merton College
\\  Oxford OX1 4JD, England} \date{} 
\maketitle
\abstract {The Lax formulation of the hyper-Hermiticity condition in
  four dimensions is used to derive a potential that generalises 
  Plebanski's second heavenly equation for
  hyper-Kahler 4-manifolds.  A class of examples of hyper-Hermitian
  metrics which depend on two arbitrary functions of two complex
  variables is given.  The twistor theory of four-dimensional
  hyper-Hermitian manifolds is formulated as a combination of the
  Nonlinear Graviton Construction with the Ward transform for
  anti-self-dual Maxwell fields.}

\noindent
\section{Complexified hyper-Hermitan manifolds}
A smooth manifold ${\cal M}$ equipped with three almost complex
structures $(I, J, K)$ satisfying the algebra of quaternions is called
hypercomplex iff the almost complex structure
\[
{\cal J}_{\l}=aI+bJ+cK
\]
is integrable for any $(a, b, c)\in S^2$. We shall use a stereographic
coordinate $\l=(a+ib)/(c-1)$ on $S^2$ which will we view as a complex
projective line $\CP^1$. 
Let $g$ be a
Riemannian metric on $\cal M$.  If $({\cal M}, {\cal J}_{\l})$ is
hypercomplex and $g({\cal J}_{\l}X, {\cal J}_{\l}Y)=g(X, Y)$ for all
vectors $X, Y$ on ${\cal M}$ then the triple $({\cal M},\; {\cal
  J}_{\l},g)$ is called a hyper-Hermitian structure.  From now on we
shall restrict ourselves to oriented four manifolds.  In four
dimensions a hyper-complex structure defines a conformal structure,
which in explicit terms is represented by a conformal frame of vector
fields $(X, IX, JX, KX)$, for any $X\in T{\cal M}$.

It is well known that this conformal structure is anti self-dual (ASD)
with the orientation determined by the complex structures.
Let $g$ be a representative of the conformal structure
defined by ${\cal J}_{\lambda}$, and let $\Sm^{A'B'}=(\Sm^{00'},
\Sm^{01'}, \Sm^{11'})$ be a basis of the space of SD two forms
${\Lambda^2}_+({\cal M})$ (see appendix for notation and conventions).
The following holds
\begin{prop}\cite{B88}
\label{boyer}
The Riemannian four manifold $({\cal M},g)$ is hyper-Hermitian
if there exists a one form $A$ (called a Lee form) depending
only on $g$ such that
\be
\label{aaa}
d\Sm^{A'B'}=-A\wedge \Sm^{A'B'}.
\ee
Moreover if $A$ is exact, then $g$ is conformally   
hyper-K\"ahler.

\end{prop}

In Section 2 we shall express  the hyper-Hermiticity  condition on the
metric in four dimensions in terms of Lax pairs of vector fields on
${\cal M}$.
The Lax formulation will be used to encode the hyper-Hermitian
geometry in a generalisation of Pleba\'nski's formalisms
\cite{Pl75}.  Some examples of hyper-Hermitian metrics are given in
Section 3.  In Section 4 we establish the twistor correspondence for
the hyper-Hermitian four-manifolds.  If ${\cal M}$ is real then the
associated twistor space is identified with a sphere bundle of
almost-complex structures and the resulting twistor theory is
well-known \cite{B88,PS93}. We will work with the complexified
correspondence and assume that ${\cal M}$ is a complex four-manifold.
The integrability conditions under which (\ref{aaa}) can hold are
$dA\in {\Lambda^2}_-({\cal M})$ so $dA$ can formally be identified
with an ASD Maxwell field on an ASD background. This will enable us to
formulate the twistor theory of hyper-Hermitian manifolds as a
non-linear graviton construction `coupled' to a Twisted Photon
Construction \cite{HW79}.
 
In Section 5 we make further remarks about the hyper-Hermitian
equation, and list some open problems. The spinor notation which
is used in the paper is summarised in the Appendix.
\section{Hyper-Hermiticity condition as an integrable system}
The hyper-Hermiticity condition on a metric $g$ can be reduced to a
system of second order PDEs for a pair of functions\footnote{ Tod has
  given a generalisation the first heavenly equation to the case of
  real hyper-Hermitian four-manifolds.}.  The Lax representation for
such an equation will be a consequence of the integrability of the
twistor distribution. We shall need the following lemma:
\begin{lemma}
Let $\nabla_{AA'}$ be four independent holomorphic vector fields on 
a four-dimensional complex manifold ${\cal M}$, and
let 
\[
L_0=\nabla_{00'}-\lambda\nabla_{01'},\;\;\;\;\;
L_1=\nabla_{10'}-\lambda\nabla_{11'},\;\;\;{\mbox{where}}\;\; 
\lambda \in \CP^1.
\]
If 
\be
\label{Lemat}
[L_0, L_1]=0
\ee
for every $\lambda$, then $\nabla_{AA'}$ is a null tetrad for a
hyper-Hermitian metric on ${\cal M}$. Every hyper-Hermitian metric
arises in this way.
\end{lemma}
{\bf Proof.}  We use the spinor notation of Penrose \& Rindler (1984).
Let $\nabla_{AA'}$ be a tetrad of holomorphic vector fields
on $\cal M$. A central result of twistor theory
\cite{Pe76, MW96} (see also Section 4 of this paper) is 
that $\nabla_{AA'}$ determines an anti self-dual conformal
structure if and only if the distribution on the primed spin bundle
$S^{A'}$ spanned by the vectors
\[
L_A=\pi^{A'}\nabla_{AA'}+\Gamma_{AA'B'C'}\pi^{A'}\pi^{B'}
\frac{\p}{\p \pi_{C'}} 
\]
is integrable. This then implies that the spin bundle is foliated by
the horizontal lifts of $\alpha$-surfaces. Here
$\pi^{A'}=\pi^{0'}o^{A'}+\pi^{1'}\iota^{A'}$ is the spinor determining
an $\alpha$-surface and is related to $\lambda=(-\pi^{1'}/\pi^{0'})$.  
From the general formula
\[
d\Sm^{A'B'}+2\Gamma_{C'}^{(A'}\wedge\Sm^{B')C'}=0,
\]
we conclude that  $\Gamma_{AA'B'C'}=-A_{A(C'}\varepsilon_{B')A'}$ for
some $A_{AA'}$ and
\[
L_A=\pi^{A'}\nabla_{AA'}+(1/2)\pi^{A'}A_{AA'}\Upsilon,
\]
where $\Upsilon=\pi^{A'}/\p \pi^{A'}$ is the Euler vector field.
We have
\begin{eqnarray}
\label{step1}
[L_A, L_B]&=&\pi^{A'}\pi^{B'}([\nabla_{AA'}, \nabla_{BB'}]
+1/2([\nabla_{BB'}, A_{AA'}\Upsilon ] -[\nabla_{AA'}, A_{BB'}\Upsilon]))\nonumber\\
&=&\pi^{A'}\pi^{B'}(
[\nabla_{AA'}, \nabla_{BB'}]+
(1/2)\varepsilon_{AB}{\nabla^C}_{(A'}A_{B')C}\Upsilon)
\nonumber\\
&=& \pi^{A'}\pi^{B'}[\nabla_{AA'}, \nabla_{BB'}]\;\;\;\;\;\;\;\;\;\;
\mbox{since $dA$ is ASD}.
\end{eqnarray}
We shall introduce the rotation coefficient $C_{ab}^c$ defined by
\[
[\nabla_a, \nabla_b]=C_{ab}^c\nabla_c,
\]
They satisfy $C_{abc}=\Gamma_{acb}-\Gamma_{bca}$. From the last
formula we can find a spinor decomposition of  $C_{abc}$, 
\[
C_{abc}=C_{ABCC'}\varepsilon_{A'B'}+C_{A'B'CC'}\varepsilon_{AB}
\]
where
\be
\label{step3}
C_{A'B'CC'}=\Gamma_{C(A'B')C'}+\varepsilon_{C'(B'}{\Gamma_{A')AC}}^{A}.
\ee
Collecting (\ref{step1}), and  (\ref{step3}) we obtain
\[
[L_A, L_B]
=\varepsilon_{AB}\pi^{A'}\pi^{B'}((1/2)A^C_{B'}{\varepsilon_{A'}}^{C'}
+{\varepsilon_{A'}}^{C'} {\Gamma_{B'D}}^{CD} )\nabla_{CC'}.
\]
We choose a spin frame $(o^{A}, \iota^A)$ constructed from two
independent solutions to the charged neutrino equation 
\[
(\nabla_{AA'}+(1/2)A_{AA'})o^{A}=(\nabla_{AA'}+(1/2)A_{AA'})\iota^{A}=0.
\]
In this frame ${\Gamma_{AA'}}^{BA}=-(1/2){A_{A'}^B}$. To obtain the equation 
(\ref{Lemat}) we project $L_A$ to the projective prime spin bundle
${\cal F}={\mathbb{P}}S_{A'}$. 
In terms of the tetrad
\be
\label{eq1}
[\nabla_{A0'}, \nabla_{B0'}]=0,
\ee
\be
\label{eq2}
[\nabla_{A0'}, \nabla_{B1'}]+[\nabla_{A1'}, \nabla_{B0'}]=0,
\ee
\be
\label{eq3}
[\nabla_{A1'}, \nabla_{B1'}]=0.
\ee
The formulation of the hyper-complex condition in formulae 
(\ref{eq1}-\ref{eq3})
was in the Riemannian case given in \cite{J95} and used in \cite{H98}.
The Lax equation 
(\ref{Lemat})
can be interpreted as the anti-self-dual Yang-Mills equations
on $\C^4$ with the gauge group $G=$Diff$({\cal M})$, reduced by four
translations in $\C^4$.
\koniec
Define $(1, 1)$ tensors ${\cal J}^{A'}_{B'}:=e^{AA'}\otimes
\nabla_{AB'}$.
As a consequence of (\ref{eq1}-\ref{eq3})
the Nijenhuis tensors 
\be
N^{A'}_{B'}(X, Y):=({\cal J}^{A'}_{B'})^2[X, Y]-
{\cal J}^{A'}_{B'}[{\cal J}^{A'}_{B'}X, Y]-
{\cal J}^{A'}_{B'}[X,{\cal J}^{A'}_{B'} Y]
+[{\cal J}^{A'}_{B'}X,{\cal J}^{A'}_{B'} Y]
\ee
vanish for arbitrary vectors $X$ and $Y$. Tensors ${\cal
J}^{B'}_{A'}$ can be treated as `complexified complex structures' on
${\cal M}$. The complex
structure ${\cal J}_{\lambda}$ on $S^{A'}$ can be conveniently
expressed as
\[
{\cal J}_{\lambda}=\pi_{A'}\tilde{\pi}^{B'}{\cal J}^{A'}_{B'},\;\;\;
\mbox{where}\;\;\;\pi_{A'}\tilde{\pi}^{A'}=1.
\]
Now we shall fix some remaining
gauge and coordinate freedom. Equations
(\ref{eq1}-\ref{eq3}) will be reduced to a coupled system of nonlinear
differential equations for a pair of functions.
\begin{prop}
\label{hyperpleban}
Let  $x^{AA'}=(x^A, w^A)$ be local null coordinates on $\cal M$ and
let $\Th^A$ be a pair of complex valued functions on $\cal M$ which
satisfy
\be
\label{hypereq}
 \frac{\p^2\Th_C}{\p x_A\p w^A}+\frac{\p\Th_B}{\p x^A}
\frac{\p^2\Th_C}{\p x_A\p x_B}=0.
\ee
Then
\be
\label{hmetric}
ds^2=dx_A\otimes dw^A+\frac{\p \Th_A}{\p x^B}dw^B\otimes dw^A
\ee
is a hyper-Hermitian metric on  $\cal M$. Conversely every 
hyper-Hermitian metric locally arises by this construction.
\end{prop}
Equation (\ref{aaa}) and its connection with 
a scalar form of (\ref{hypereq}) was investigated 
by different methods in \cite{FP76} in the context of weak heavenly
spaces. 
Other integrable equations associated to 
Hyper-Hermitian manifolds have been studied 
in \cite{GS98}.\\
{\bf Proof.}
Choose a conformal factor such that
$A_{AA'}=o_{A'}A_A$ for some $o_{A'}$ and  $A_A$.
This can be done since the two form $\Sm^{1'1'}$ is simple and therefore
equation (\ref{aaa}) together with the Frobenius theorem imply the existence of
the conformal factor such that $d\Sm^{1'1'}=0$. Hence,
using the Darboux's theorem,
one
can introduce canonical coordinates $w^A$ such that
\[
\Sm^{1'1'}=(1/2)\varepsilon_{AB}dw^{A}\wedge dw^B,
\]
and choose un unprimed spin frame so that $o_{A'}e^{AA'}=dw^A$.
Coordinates $w^A$  parametrise the space of null surfaces tangent to 
$o^{A'}$, i.e. $o^{A'}\nabla_{AA'}w^B=0$.
Consider  
\[
{\cal J}^{1'}_{0'}=o^{B'}dw^A\otimes \nabla_{AB'}
\]
The tensor
${\cal J}^{1'}_{0'}$ is a degenerate complex structure. Therefore 
$({\cal J}^{1'}_{0'})^2=0$ where ${\cal J}^{1'}_{0'}$ is now
thought of as a differential operator acting on forms. Let $h$ be
a function on ${\cal M}$. Then
\[
{\cal J}^{1'}_{0'}\hook d ({\cal J}^{1'}_{0'}(dh))=0\;\;\;\;
\mbox{implies that}\;\;
[\nabla_{A0'}, \nabla_{B0'}]=0,
\]
and our choice of the spin frame is consistent with
(\ref{eq1}-\ref{eq3}). By applying the Frobenius theorem we can find
coordinates $x^A$ such that 
\[
\label{mt}
\nabla_{A0'}=\frac{\p}{\p x^A},\;\;\;\;\;\;
\nabla_{A1'}=
\frac{\p}{\p w^A}-{\Th_A}^B\frac{\p}{\p x^B}.
\] 
Using  equation (\ref{eq2}), 
we deduce the existence of a 
potential
$\Th_A$ such that
${\Th_A}^B=\nabla_{A0'}\Th^B$.
Now (\ref{eq3}) gives the
field equations (\ref{hypereq})
\[
\frac{\p^2\Th_C}{\p x_A\p w^A}+\frac{\p\Th_B}{\p x^A}
\frac{\p^2\Th_C}{\p x_A\p x_B}=0.
\]
The dual frame  is
\[
e^{A0'}=dx^A+\frac{\p \Th^A}{\p x^B}dw^B,\;\;\;\;\;\;e^{A1'}=dw^A,
\]
which justifies formula (\ref{hmetric}).\koniec
In the adopted gauge, the Maxwell potential is 
\[
A=\frac{\p^2 \Th^B}{\p x^A \p x^B}dw^A
\]
and $\nabla^aA_a=0$ i.e. this is a `Gauduchon gauge'.
Electromagnetic gauge transformations on $A$ correspond to conformal
rescalings of the tetrad (which preserve the hypercomplex structure).
The second heavenly equation (and therefore the hyper-K\"ahler condition)
follows from (\ref{hypereq}) if in addition $\nabla_{A0'}\Th^A=0$.
This  condition  guarantees the existence of a
scalar function $\Th$ which satisfies the second Pleba{\'n}ski equation
\[
\frac{\p^2\Th}{\p w^A\p x_A}+
\frac{1}{2}\frac{\p^2\Th}{\p x^B\p x^A}\frac{\p^2\Th}{\p x_B\p x_A}=0
\]
such that $\Th^A={\nabla^A}_{0'}\Th$.
In this case  $A$ is  exact so can be gauged away by a conformal
rescaling.
\section{Examples}
We look for solutions to (\ref{hypereq}) for which the linear and
nonlinear terms vanish separately, ie.
\be
\label{specialcase}
\frac{\p^2\Th_C}{\p x_A\p w^A}=
\frac{\p\Th_B}{\p x^A}
\frac{\p^2\Th_C}{\p x_A\p x_B}=0.
\ee
Put $w^A=(w, z), x^A=(y,-x)$. 
A simple class of  solutions to (\ref{specialcase}) is provided by
\[
\Th_0=ax^l,\;\;\Th_1=by^k,\;\;\;\;\;\;\;\;\;\;\;\;k,l\in\Z,\;\;a,b\in\C.
\]
The corresponding metric and the Lee form are
\be
\label{ex2}
ds^2=dw\otimes dx+dz\otimes dy+ (alx^{l-1}+bky^{k-1})dw\otimes dz,
\ee
\[
A= {b\left (k-1\right )k{y}^{k-2}dw}- {a\left (l-1\right )l{x}^{l-1}dz}.
\]
From calculating the invariant
\[
C_{ABCD}C^{ABCD}= (3/2)abk(k-1)(k-2)l(l-1)(l-2) x^{l-3}y^{k-3}
\]
we conclude that the metric (\ref{ex2}) is in general of type $I$ or
$D$ 
(or type $III$ or$N$ if $a$ or $b$ vanish, or $k<3$ or $l<3$).
\subsection{Hyper-Hermitian elementary states}
A  more interesting class of solutions
(which  generalise the metric of Sparling and Tod
\cite{ST79} to the hyper-Hermitian case) is given by
\be
\label{ex}
\Th_C=\frac{1}{x_Aw^A}F_C(W^A),  
\ee
where $W^A=w^A/(x_Bw^B)$ and $F_C$ are two  arbitrary
complex functions of two complex variables. The corresponding metric is
\[
ds^2=dx_A\otimes dw^A+\frac{1}{(x_Aw^A)^2}\Big(F_C+\frac{w^B}{(x_Aw^A)}
\frac{\p F_C}{\p W^B}\Big)dw^C\otimes(w_Adw^A).
\]
This metric is singular at the light-cone of the origin. The
singularity may be moved to infinity if we introduce new
coordinates $X^A=x^A/(x_Bw^B), W^A=w^A/(x_Bw^B)$ and
rescale the metric  by $(X_AW^A)^2$.
This yields
\be
\hat{{ds}^2}=dX_A\otimes dW^A+\Big(F_B+W^C\frac{\p F_B}{\p W^C}\Big)
\Big((X_AW^A)dW^B-W^Bd(X_AW^A)\Big)\otimes(W_AdW^A)
\ee
and
\[
A=-\Big(3W^AF_A+5W^AW^B\frac{\p F_A}{\p W^B}+W^AW^BW^C\frac{\p^2
F_A}{\p W^B\p W^C}\Big)W_DdW^D.
\]
The metric of Sparling and Tod corresponds to setting $F_A=W_A$.

 Let us consider the particular case $F_A=(aW^kZ^l, bW^mZ^n)$. The metric is
\be
\label{exa3}
ds^2=dw\otimes dx+dz\otimes dy+\Big( \frac{ a(k+l+1)w^kz^l}{(wx+zy)^{k+l+2}}dw+
\frac{b(m+n+1)w^mz^n}{(wx+zy)^{m+n+2}}dz \Big)\otimes(wdz-zdw).
\ee
If $ a=-b, l=n+1, k=m-1$ then $\Th_A=\nabla_{A0'}\Th$  where
$\Th=-aw^kz^{l-1}(wx+zy)^{-(k+l)}$. 
For these values of parameters the metric is hyper-K\"ahler and of
type $N$. 
 
 Some solutions to (\ref{specialcase}) have real Euclidean slices. 
For example 
\[
\Th_0=-{\frac {y\left (2\,wx+zy\right )}{{w}^{2}\left (wx+zy\right
)^{2}}},\;\;\;
\Th_1=-{\frac {{y}^{2}}{w\left (wx+zy\right )^{2}}}
\]
with $w=\bar{x}, z=\bar{y}$
yield a solution of type $D$, which is conformal to the Eguchi-Hanson
metric.
\section{The twistor construction}
In this section we shall establish the  correspondence between
complexified hyper-Hermitian four manifolds and three dimensional
twistor spaces with additional structures. We shall also look at
examples given in Section 3 from  the twistor point of view.
 We begin with recalling
basic facts about the twistor  correspondences for 
ASD spaces \cite{Pe76, MW96}.

Define $\a$-surfaces as null self-dual two-dimensional surfaces in
$\cal M$. The correspondence space $\cal F$ is a set of pairs $(x,
\lambda)$ where $x\in \cal M$ and $\l \in \CP^1$ parametrises
$\a$-surfaces through $x$ in $\cal M$.  We represent $\cal F$ as the
quotient of the primed-spin bundle $S^{A'}$ with fibre coordinates
$\pi_{A'}$ by the Euler vector field $\pi^{A'}/\p \pi^{A'}$ so that
the fibre coordinates are related by $\lambda=\pi_{0'}/\pi_{1'}$.  The
space $\cal F$ possesses a natural two dimensional distribution
(called the twistor distribution, or the Lax pair, to emphasise the
analogy with integrable systems).  The Lax pair on ${\cal F}$ arises
as the image under the projection $TS^{A'}\longrightarrow T{\cal F}$
of the distribution spanned by
\[
\pi^{A'}\nabla_{AA'}+\Gamma_{AA'B'C'}\pi^{A'}\pi^{B'}\frac{\p}{\p
\pi_{C'}}
\]
and is given by
\be
\label{laxpair}
L_A=(\pi_{1'}^{-1})(\pi^{A'}\nabla_{AA'}+f_A\p_{\lambda}),\;\;\;\;\mbox{
where }\;\; f_A=(\pi_{1'}^{-2})\Gamma_{AA'B'C'}\pi^{A'}\pi^{B'}\pi^{C'}.
\ee
The integrability of the twistor
distribution is equivalent to $C_{A'B'C'D'}=0$, the vanishing of 
the self-dual Weyl spinor. The  twistor space arises as 
a factor space of $\cal F$ by the twistor distribution.
This leads to a double fibration
\be
\label{doublefib}
{\cal M}\stackrel{p}\longleftarrow 
{\cal F}\stackrel{q}\longrightarrow {\cal PT}.
\ee

The existence of $L_A$ can also be deduced directly from the
correspondence with $\cal PT$. The basic twistor correspondence
\cite{Pe76} states that points in ${\cal M}$ correspond in ${\cal PT}$
to rational curves with normal bundle ${\cal O}^{A}(1)={\cal
  O}(1)\oplus {\cal O}(1)$.  Let $l_x$ be the line in $\cal PT$ that
corresponds to $x\in {\cal M}$. The normal bundle to $l_x$ consists of
vectors tangent to $x$ (horizontally lifted to $T_{(x,\lambda)}{\cal
  F}$) modulo the twistor distribution. Therefore we have a sequence
of sheaves over $\C P^1$
\[
0\longrightarrow D \longrightarrow \C^4 \longrightarrow
{\cal O}^A(1)\longrightarrow 0.
\]
The map $\C^4 \longrightarrow {\cal O}^A(1)$ is given by
$V^{AA'}\longrightarrow V^{AA'}\pi_{A'}$.  Its kernel consists of
vectors of the form $\pi^{A'}\lambda^A$ with $\lambda^A$ varying. The
twistor distribution is therefore $D=O(-1)\otimes S^{A}$ and so $L_A$,
the global section of $\Gamma(D\otimes {\cal O}(1)\otimes S_{A})$, has
the form (\ref{laxpair}).


 We have 
\begin{prop}
\label{twmaci}
Let $\cal PT$ be a three-dimensional complex manifold with the following 
structures
\begin{itemize}
\item[(A)] a projection $\mu :{\cal PT}\longrightarrow \CP^1$,
\item[(B)] a four complex dimensional family of sections with
a normal bundle ${\cal O}(1)\oplus {\cal O}(1)$.
\end{itemize}
Then the moduli space ${\cal M}$ of sections of $\mu$ is equipped with 
hyper-Hermitian structure. Conversely given a hyper-Hermitian
four-manifold there will always exists a corresponding twistor space 
satisfying conditions $(A)$ and $(B)$.
\end{prop}
\noindent
{\bf Remarks}
\begin{itemize}
\item [(i)] 
The integrability conditions under which (\ref{aaa}) can hold  are  
$dA\in {\Lambda^2}_-({\cal M})$
so $dA$ can  be identified with an ASD Maxwell field on an ASD
background.
This suggests that hyper-Hermitian manifolds can be studied with
respect to  the Twisted Photon Construction \cite{HW79}, associated
with $dA$.
Let $K=\Lambda^3({\cal PT})$ be the canonical line bundle.
Proposition 4 is different from the original Nonlinear Graviton
construction because the line bundle $L:=K^*\otimes{\cal O}(-4)$,
where ${\cal O}(-4)=\mu^*(T^*\CP^1)^2$,
is in general nontrivial over ${\cal PT}$. It is the twisted photon
line bundle associated with $dA$. 
\item[(ii)]
If ${\cal M}$ is compact then it follows from Hodge theory that $dA=0$
and the hyper-Hermitan structure
is locally conformally hyper-K\"ahler. We focus on the non-compact case.
\item[(iii)] If ${\cal M}$ is real then ${\cal PT}$ is equipped
with an antiholomorphic 
involution preserving $(A)$ and we recover a result 
closely related to one
of Petersen and
Swann \cite{PS93} who constructed a twistor space corresponding to a 
real four-dimensional ASD Einstein--Weyl metric with vanishing scalar
curvature. 
\item[(iv)]
The correspondence is preserved under
holomorphic deformations of ${\cal PT}$ which preserve $(A)$.
\end{itemize}

\noindent
{\bf Proof .} Let
$\l=\pi_{0'}/\pi_{1'}$ be an affine coordinate on $\CP^1$.
$\cal PT$ can be covered by two sets, $U$ and $\tilde{U}$ with
$|\lambda|< 1+\epsilon$ on $U$ and $|\lambda|>1-\epsilon$ on
$\tilde{U}$ with $(\om^A,\lambda)$ being the coordinates on $U$ and
$(\tilde{\om}^A,\lambda^{-1})$ on $\tilde{U}$.  The twistor space
$\cal PT$ is then
determined by the transition function ${\tilde{\om}}^B=
{\tilde{\om}}^B(\om^A, \pi_{A'})$ on $U\cap \tilde{U}$.
Let $l_x$ be the line in $\cal PT$ that
corresponds to  $x\in {\cal M}$ and let $Z\in\cal PT$ lie on $l_x$. 
We denote by $\cal F$ the correspondence space
${\cal PT}\times {\cal M}|_{Z\in l_x}= {\cal M}\times\CP^1$ and
 use the double fibration picture (\ref{doublefib}).
 
Consider the line bundle 
\[
L= K^*\otimes{\cal O}(-4)
\] 
over ${\cal PT}$ given by the transition function
$f=\mbox{det}(\p\tilde{\om}^A/\p \om^B)$. When pulled back to $\cal F$
it satisfies
\[
L_Af=0.
\]
Since $H^1({\cal F}, {\cal O})=0$, we can
perform the splitting $f=h_0h_{\infty}^{-1}$.
By the standard Liouville arguments (see \cite{HW79}) we deduce that
\be
\label{max}
h_0^{-1}L_A(h_0)=h_{\infty}^{-1}L_A (h_{\infty})=-(1/2)A_A
\ee
where $A_A=A_{AB'}\pi^{B'}$ is global on ${\cal F}$.
The integrability conditions imply that
$F_{AB}=\nabla_{A'(A}A_{B)}^{A'}$ is an ASD Maxwell field on the ASD
background. The one-form  $A=A_{AA'}e^{AA'}$ is a Maxwell potential.
 The canonical line bundle of ${\cal PT}$ is
$K={\cal O}(-4)\otimes L^*$. 
To obtain a global, line bundle valued three-form on ${\cal PT}$ one 
must tensor the last equation
with ${\cal O}(4)\otimes L$. We pick a global section 
$\xi\in\Gamma(K\otimes{\cal O}(4)\otimes L)$ and restrict $\xi$ to $l$ 
\be
\label{volume}
\xi|_l=\Sm_{\l}\wedge \pi_{A'}d\pi^{A'}
\ee
where $\pi_{A'}d\pi^{A'}\in \Om^1\otimes{\cal O}(2) $. A two-form 
\be
\label{stform}
\Sm_{\l} \in \Gamma(\Lambda^2(\mu^{-1}(\lambda))\otimes {\cal O}(2)\otimes L)
\ee
is defined on vectors vertical with respect to $\mu$ 
by $\Sm_{\l}(U, V)\pi_{A'}d\pi^{A'}=\xi(U, V, ...)$.
Let $p^*\Sm_{\l}$ be the pullback of
$\Sm_{\l}$ to ${\cal F}$. Note that if
\[ 
A\longrightarrow A-d\phi\;\; \mbox{(gauge transformation on $L$)}
\;\;\;\;\;\;\; \mbox{then}\;\;\;\;\;\; 
p^*\Sm_{\l}\longrightarrow e^{\phi}p^*\Sm_{\l}.
\]
Let $p^*\Sm_{\l}$ be defined over $U$ and  $p^*\widetilde\Sm_{\l}$ over 
$\tilde U$. We have $f(p^*\Sm_{\l})=p^*\widetilde\Sm_{\l}$. By  definition,
$p^*\Sm_{\l}$ descends to the twistor space, i.e.,
\be
\label{proj}
{\cal L}_{L_A}(p^*\Sm_{\l})=0.
\ee
We make use of the splitting formula, and define (on ${\cal F}$)
$\Sm_0=h_0(p^*\Sm_{\l})$.
The line bundle valued two-form $\Sm_0$ is a globally defined object on
$\cal F$, and therefore it is equal to $\pi_{A'}\pi_{B'}\Sm^{A'B'}$.
Note that $\Sm_0$ does not descend to ${\cal PT}$.
Fix $\l\in \CP^1$ (which gives 
a copy ${\cal M}_{\l}$ of ${\cal M}$ in ${\cal F}$) and apply
(\ref{proj}). This yields
\[
{\cal L}_{L_A}\Sm_0=h_0^{-1}{L_A}(h_0)\Sm_0.
\]
After some work we obtain formula (\ref{aaa}): 
\be
\label{hypcom}
d\Sm^{A'B'}=-A\wedge \Sm^{A'B'}.
\ee
The integrability conditions for the last equation are guaranteed by the
existence of solutions to (\ref{max}).
Equation (\ref{hypcom})  and  the forward part of Proposition \ref{boyer} 
imply that $\cal M$ is equipped with 
hyper-Hermitan structure. If the line bundle $L$ over $\cal PT$ is
trivial, then $\cal M$  is  conformally hyper-K\"ahler.

 Now we discuss the converse problem of recovering various structures
on ${\cal PT}$   from the geometry of $\cal M$.
Let $\cal M$ by  a hyper-Hermitian four-manifold. Therefore
$C_{A'B'C'D'}=0$ and there exists a twistor space satisfying 
Condition $(A)$.
Equation (\ref{hypcom}) implies that  $F=dA$ is an ASD Maxwell field,
and we can solve
\[
\pi^{A'}(\nabla_{AA'}+(1/2)A_{AA'})\rho=0
\]
on each $\a$-surface (self-dual, two dimensional null surface in
${\cal M}$). We define fibres of $L$ as one-dimensional spaces of
solutions to the last equation. The solutions on $\a$-surfaces
intersecting at $p\in \cal M$ can be compared at one point, so $L$
restricted to a line $l_x$ in ${\cal PT}$ is trivial.  In order to
prove that ${\cal PT}$ is fibred over $\CP^1$ notice that equation
$\pi^{A'}(\nabla_{AA'}+(1/2)A_{AA'})\pi_{B'}=0$ implies
$\pi^{A'}\nabla_{AA'}\lambda =0$, so $\lambda$ and $1/\lambda$ descend
to give meromorphic functions on twistor space and defines the map
${\cal PT}\rightarrow\CP^1$.
\koniec
\subsection{Examples}
In this Subsection we shall give the twistor correspondence for
the family of hyper-Hermitian metrics (\ref{exa3}) found in Section 3.
First we shall look at the passive twistor constructions of $\Th_C$ 
by the contour integral formulae. It will turn out that 
$\Th_C$ are examples of  Penrose's elementary states. 
Then we explain how the cohomology classes corresponding to $\Th_C$
can be used to deform a patching description of ${\cal PT}$. The deformed
twistor space will, by proposition 4,
give rise to the metric (\ref{exa3}). Both passive and active
constructions in this subsection use methods developed by Sparling in
his twistorial treatment of the Sparling-Tod metric.
 
Parametrise a section of $\mu:{\cal PT}\longrightarrow \CP^1$
by the coordinates
\[
x^{AA'}:=\frac{\p \om^A}{\p \pi_{A'}}{\Big |}_{\pi_{A'}=o_{A'}}=
\left (
\begin{array}{cc}
y&w\\
-x&z
\end{array}
\right ),\;\;\; 
{\mbox{so that}}\; x^{A1'}=w^A=(w, z),\;x^{A0'}=x^A=(y,-x).
\] 

Let us consider the particular case $F_A=(aW^kZ^l, bW^mZ^n)$
discussed in Subsection 3.1.
We work on the non-deformed  twistor space $\cal PT$ with homogeneous
coordinates $(\om^A, \pi_{A'})$. On the primed spin bundle
$\om^0=\pi_{1'}(w+\lambda y), \om^1=\pi_{1'}(z-\lambda x)$. 
Consider two twistor functions (sections of $H^1(\CP^1, {\cal{O}}(-2)$) 
\[
h_0=(-1)^k a\frac{(\pi_{0'})^{k+l}}{(\om^0)^{l+1}(\om^1)^{k+1}},\;\;\;
h_1=(-1)^m b\frac{(\pi_{0'})^{m+n}}{(\om^0)^{n+1}(\om^1)^{m+1}}.
\]
where $a, b\in \C$ and $k, l, m, n \in \Z$ are constant  parameters. 
Then
\[
\Th_A(w,z,x,y)=\frac{1}{2\pi i} \oint_{\Gamma}h_A(\om^{B}, \pi_{B'})
\pi_{A'}d\pi^{A'}.\]
Here  $\Gamma$  is a contour in $l_x$, the $\CP^1$  that corresponds
to $(w, z, x, y)\in {\cal M}$. It
separates the two poles of the integrand. To find
$\Th^A$ we compute the residue at one of these poles, which gives
\be
\label{elstates}
\Th_0=a\frac{w^kz^l}{(wx+zy)^{k+l+1}},\;\;\;
\Th_1=b\frac{w^mz^n}{(wx+zy)^{m+n+1}}, 
\ee
and hence the metric (\ref{exa3}).   

Now we shall use $h_A$ to deform the complex structure of 
${\cal PT}$.
We change the standard patching relations by setting
\[
\tilde{\om}^A=f^A(\om^{A},t)
\]
where 
$t$ is a deformation parameter and $f^A$ is determined by the
deformation equations
\be
\label{defoeq}
\frac{df^0}{dt}=\frac{b\pi_{0'}^{m+n+3}}{(\tom^0)^{n+1}(\tom^1)^{m+1}}
(-1)^m,\;\;
\frac{df^1}{dt}=\frac{a\pi_{0'}^{k+l+3}}{(\tom^0)^{l+1}(\tom^1)^{k+1}}
(-1)^{k+1}.  
\ee 
This equation has first integrals.  If $ a=-b, l=n+1,
k=m-1$ then (\ref{defoeq}) imply that $\om^0\om^1=\tom^0\tom^1$ is a
global twistor function.  When pulled back to the spin bundle this can
be expressed as $P_{A'B'}\pi^{A'}\pi^{B'}$, and the corresponding
metric admits a null Killing vector $K_{AA'}$ given by
\[
\nabla_{AC'}P_{A'B'}=K_{A(A'}\varepsilon_{B')C'}.
\]
Assume that $n+1\neq l$, and   $k+1\neq m$. Then a
first integral of (\ref{defoeq})
\[
Q=\frac{a(\pi^{0'})^{k+l+3}(-1)^{k+1}}{n+1-l}(\om^0)^{n+1-l}
+\frac{b(\pi^{0'})^{m+n+3}(-1)^{m+1}}{k+1-m}(\om^1)^{k+1-m}
\]
is given by a function homogeneous of degree $k+n+4$. 
Its pull backs to ${\cal F}$  (which we also denote $Q$) satisfies
$L_A(Q)=0$. This implies the existence of a Killing spinor of valence 
$(0,k+n+4)$ on ${\cal M}$.
\section{Further Remarks}
\subsection{Symmetries}
The equation (\ref{hypereq}) has the obvious first
integral given by functions $\Lambda_C$ which satisfy
\[ 
\frac{\p\Th_C}{\p w^A}+\frac{\p\Th_B}{\p x^A}
\frac{\p\Th_C}{\p x^B}=\frac{\p \Lambda_C}{\p  x^A}.
\]
It is implicit from the twistor construction that equation
(\ref{hypereq}) has infinitely many first integrals given by hidden
symmetries. These will be studied (and the associated hierarchy of
equations \cite{DM96}) in a subsequent paper.  Here we give a
description of those symmetries that correspond to the pure gauge
transformations.

Let $M$ be a vector field on $\cal M$. Define
$\delta_M^0\nabla_{AA'}:=[M, \nabla_{AA'}]$.  This is a pure gauge
transformation corresponding to the addition of ${\cal L}_Mg$ to the
space-time metric.

Once a coordinate system 
leading to equation (\ref{hypereq})
has been selected, the  field equations will
not be invariant under all the diff$(\cal M)$ transformations.
We restrict ourselves to transformations that preserve 
the canonical structures on ${\cal M}$, namely
\[
\Sm^{1'1'}=(1/2)dw_A\wedge dw^A,\;\;\; \mbox{and}\;\;\;\
{\cal J}^{1'}_{0'}=dw^A\otimes\frac{\p}{\p x^A}. 
\]
The
condition ${\cal L}_M\Sm^{0'0'}={\cal L}_M{\cal J}^{1'}_{0'} =0$ implies that
$M$ is given by
\[
 M= \frac{\p h}{\p w_A }  \frac{\p}{\p w^A} +\Big( g^A-
x^B\frac{\p^2 h}{\p w_A \p w^B} \Big)\frac{\p}{\p x^A}
\]
where $h=h(w^A)$ and $g^A=g^A(w^B)$.
Space-time is now viewed as a tangent bundle ${\cal M}=T{\cal
N}^2$ 
with $w^A$ being coordinates on the 
two-dimensional complex manifold ${\cal N}^2$. The full diff$({\cal M})$ 
symmetry breaks down to  sdiff$({\cal N}^2)$  
which acts on ${\cal M}$ by Lie lift. Let $\delta_M^0\Th$
corresponds to $\delta_M^0\nabla_{AA'}$ by
\[
\delta_M^0\nabla_{A1'}=\frac{\p\delta_M^0\Th^B}{\p x^A}\frac{\p}{\p x^B}.
\]
The `pure gauge' elements are 
\[
\delta^0_M \Th^B={\cal L}_M(\Th^B)+F^B-x^A\frac{\p g^B}{\p w^A}
+x^Ax^C \frac{\p^2 h}{\p w^A \p w^C \p
w_B}
\]
where $F^B, g^A$ and  $h$ are functions of $w^B$ only.
\subsection{$gl(2, \C)$ connection}
A natural connection which arises in hyper-Hermitian geometry
is the Obata connection \cite{O56}.
In this subsection we discuss other possible choices of connections
associated with hyper-Hermitian geometry. 
 We shall motivate our choices by considering the conformal rescalings
of the null tetrad.
  The first Cartan structure equations are
\[
de^{AA'}=e^{BA'}\wedge{\Gamma^{A}}_{B}+e^{AB'}\wedge{\Gamma^{A'}}_{B'}.
\]
Rescaling $e^{AA'}\longrightarrow \hat{e}^{AA'}= e^{\phi}e^{AA'}$ yields
\[
d\hat{e}^{AA'}=\hat{e}^{BA'}\wedge{\Gamma^{A}}_{B}+\hat{e}^{AB'}\wedge{\Gamma^{A'}}_{B'}+d\phi\wedge\hat{e}^{AA'}.
\]
The last equation can be interpreted in (at least) three different ways;
\begin{itemize}
\item[(a)]
Introduce the torsion three-form by 
$T=\ast(d\phi)=T_{abc}\hat{e}^a\wedge\hat{e}^b \wedge\hat{e}^c$. Then
\[
d\hat{e}^a+{\Gamma^a}_b\wedge \hat{e}^b=T^a
\]
where $T^a=(1/2)T^a_{bc}\hat{e}^b \wedge\hat{e}^c$.
\item[(b)]
Use the torsion-free $sl(2, \C)\oplus\widetilde{sl}(2, \C)$ spin connection
\[
\Gamma_{AB}\longrightarrow \Gamma_{AB}+1/4\ast(d\phi\wedge\Sm_{AB}),\;\;\;
\Gamma_{A'B'}\longrightarrow \Gamma_{A'B'}+
1/4\ast(d\phi\wedge\Sm_{A'B'}),\;\;\;
\]
\item[(c)]
Work with the torsion-free $gl(2, \C)\oplus \widetilde{gl}(2, \C)$ connection
\[
{\G}_{AB}=\Gamma_{AB}+a\varepsilon_{AB}d\phi,\;\;\;
{\G}_{A'B'}=\Gamma_{A'B'}+ (1-a)\varepsilon_{A'B'}d\phi
\]
with $\Gamma_{AB}=\Gamma_{(AB)}\in {sl}(2,\C)\otimes
\Lambda^{1}(T^*{\cal M})$,  
$\Gamma_{A'B'}=\Gamma_{(A'B')}\in \widetilde{sl}(2,\C)\otimes
\Lambda^{1}(T^*{\cal M})$ and $a\in \C$.
This leads to
\[
d\hat{e}^a+{\G^a}_b\wedge \hat{e}^b=0
\] 
where ${\G}_{ab}=\Gamma_{ab}+\varepsilon_{A'B'}\varepsilon_{AB}d\phi$.
The structure group reduces to
\[
sl(2, \C)\oplus \widetilde{sl}(2,
\C)\oplus u(1)\subset
gl(2, \C)\oplus \widetilde{gl}(2, \C).
\]
\end{itemize}
For (complexified) hyper-Hermitian four-manifolds
$d\phi$ is replaced by the Lee form $-A$ in the above formulae.  
The possibility $(a)$ would then correspond to the heterotic
geometries studied by physicists in connection with $(4, 0)$
supersymmetric $\sigma$-models (see \cite{CTV96} and references
therein). Choice $(b)$ is what we have
used in this paper. Let us make a few remarks about the possibility
$(c)$.
 
 Equation (\ref{aaa}) implies that $a=1/2$ and
\[
\G_{AB}=\Gamma_{AB}-1/2\varepsilon_{AB}A,\;\;\;
\G_{A'B'}=-1/2\varepsilon_{A'B'}A
\]
with $\Gamma_{AB}=\Gamma_{(AB)}\in  sl(2,\C)$. In the adopted
coordinate system
\[
\Gamma_{AA'BC}=-o_{A'}\Big(\nabla_{(A0''}\nabla_{B0'}\Th_{C)}+
\frac{1}{2}\varepsilon_{BC}\frac{\p \Th^D}{\p x^A \p x^D}\Big),\;\;\;\;
\Gamma_{AA'B'C'}=
-\frac{1}{2}o_{A'}\varepsilon_{B'C'}\frac{\p \Th^D}{\p x^A \p x^D}.
\]
The curvatures of $\G_{AB}$ and  $\G_{A'B'}$ are 
\[
{\R^A}_B=d{\G^A}_B+{\G^A}_C\wedge{\G^C}_B=
{R^A}_B-{1}/{2}{\varepsilon^A}_BF,\;\;\;\;
{\R^{A'}}_{B'}=-{1}/{2}{\varepsilon^{A'}}_{B'}F
\]
where $F=dA$ is an ASD two form.
It would be interesting to investigate this possibility
with connection to $gl(2, \C)$ formulation of Einstein-Maxwell
equations \cite{Pl74}, and its Lagrangian description \cite{R94}. 
\subsection{Reductions}
Hyper-Hermitian four-manifolds which admit a tri-holomorphic vector
field were recently studied in \cite{CTV96} and \cite{GT98}. It would
be interesting to look at the case of a general Killing vector taking
the equation (\ref{hypereq}) as a starting point.  One might also
consider reduction of real slices with $(++--)$ signature to obtain an
`evolution' form of Einstein-Weyl equations for metrics of signature
$(+--)$.
\section{Acknowledgement}
I am grateful to Dr Lionel Mason and Dr Paul Tod for helpful
discussions, and to Merton College for a Palmer Senior Scholarship.
\section{Appendix}
We shall use the conventions of
Penrose and Rindler \cite{PR86}:
$a,b,...$ are four-dimensional space-time indices and $A, B, ..., A', B',
...$ are two-dimensional spinor indices. The tangent space  at each
point of ${\cal M}$ is isomorphic to a tensor product of two spin spaces
\[
T^a{\cal M}=S^A\otimes S^{A'}.
\]
Spin dyads $(o^A, \iota^{A})$ and 
$(o^{A'}, \iota^{A'})$ span $S^A$ and $S^{A'}$ respectively.
The spin spaces $S^A$ and $S^{A'}$ are equipped with symplectic forms
$\varepsilon_{AB}$
and $\varepsilon_{A'B'}$ such that
$\varepsilon_{01}=\varepsilon_{0'1'}=1$.
These anti-symmetric objects are used to  raise and lower 
the spinor indices. We shall use the normalised spin frames, which
implies that
\[
o^B\iota^C-\iota^Bo^C=\varepsilon^{BC},\;\;\;
o^{B'}\iota^{C'}-\iota^{B'}o^{C'}=\varepsilon^{B'C'}.
\] 
 Let $e^{AA'}$ be the null tetrad of one forms on ${\cal M}$ and let
$\nabla_{AA'}$ be the frame of dual vector fields. 
The orientation is fixed by setting 
$\nu=e^{01'}\wedge e^{10'}\wedge e^{11'}\wedge e^{00'}$.
The local basis $\Sm^{AB}$ and $\Sm^{A'B'}$ of  spaces of ASD and SD
two forms are defined by
\[
e^{AA'}\wedge e^{BB'}=\varepsilon^{AB}\Sm^{A'B'}+\varepsilon^{A'B'}\Sm^{AB}. 
\]
The Weyl tensor decomposes into ASD and SD part 
\[
C_{abcd}=\varepsilon_{A'B'}\varepsilon_{C'D'}C_{ABCD}+
\varepsilon_{AB}\varepsilon_{CD}C_{A'B'C'D'}.
\]
The first Cartan structure equations are 
\[
de^{AA'}=e^{BA'}\wedge{\Gamma^{A}}_{B}+e^{AB'}\wedge{\Gamma^{A'}}_{B'},
\]
where $\Gamma_{AB}$ and $\Gamma_{A'B'}$ are the  $SL(2, \C)$
and $\widetilde{SL}(2, \C)$
spin connection one forms symmetric in their indices, and
\[
 \Gamma_{AB}=
\Gamma_{CC'AB}e^{CC'},\;\;\Gamma_{A'B'}=\Gamma_{CC'A'B'}e^{CC'}
,\;\;\; \Gamma_{CC'A'B'}=o_{A'}\nabla_{CC'}\iota_{B'}-
\iota_{A'}\nabla_{CC'}o_{B'}.
\]
The curvature of the spin connection
\[
{R^A}_B=d{\Gamma^A}_B+{\Gamma^A}_C\wedge{\Gamma^C}_B
\]
decomposes as
\[
{R^A}_B={C^A}_{BCD}\Sm^{CD}+(1/12)R{\Sm^{A}}_{B}+{\Phi^A}_{BC'D'}\Sm^{C'D'}
\]
and similarly  for ${R^{A'}}_{B'}$. Here $R$ is the Ricci scalar and 
$\Phi_{ABA'B'}$ is the trace-free part of the Ricci tensor.

For convenience we express various spinor objects on $\cal M$ in terms
of $\Th_A$.
\begin{eqnarray*}
\mbox{Tetrad}& &e^{A0'}=dx^A+\frac{\p \Th^A}{\p x^B}dw^B,\;\;e^{A1'}=dw^A,\\
\mbox{dual tetrad}&&
\nabla_{A0'}=\frac{\p}{\p x^A},\;\;\;\;
\nabla_{A1'}=\frac{\p}{\p w^A}-\frac{\p \Th^B}{\p x^A}\frac{\p}{\p x^B},\\
\mbox{metric determinant}&&det (g)=1\\
\mbox{Weyl spinors}&&C_{A'B'D'E'}=0,\;\;
C_{ABCD}=\nabla_{(A0'}\nabla_{B0'}\nabla_{C0'}\Th_{D)},\\
\mbox{spin connections}&&\Gamma_{AA'BC}=-\frac{1}{2}
o_{A'}(\nabla_{(B0'}\nabla_{C0'}\Th_{A)}+\nabla_{B0'}\nabla_{C0'}\Th_A),\\
&&\Gamma_{AA'B'C'}=-\frac{\p^2 \Th^B}{\p x^B\p x^A}o_{(B'}\varepsilon_{C')A'},\\
\mbox{Lee form}&& A=\frac{\p^2 \Th^B}{\p x^B\p x^A}dw^A, \\
\mbox{wave operator}&&
\square_g=A^{a}\p_{a}+\nabla^{A}_{1'}\nabla_{A0'}=
\frac{\p^2}{\p x_A\p w^A}+
\frac{\p^2\Th_B}{\p x_A\p x_B}\frac{\p}{\p x^A}
+\frac{\p\Th^A}{\p x_B}\frac{\p}{\p x^A}\frac{\p}{\p x^B},\\
\mbox{ Ricci scalar}&& R=1/12(\nabla^aA_a+A_aA^a)=0. 
\end{eqnarray*}
The last formula follows because $A$ is null and satisfies the Gauduchon gauge.

\end{document}